\documentclass[12pt]{amsart}
\usepackage{amsthm,amssymb,latexsym}
\usepackage{graphics, epsfig}

\theoremstyle{plain}
\newtheorem{thm}{Theorem}

\newtheorem{lem}{Lemma}
\newtheorem{cor}{Corollary}

\newtheorem{conjecture}{Conjecture}

\theoremstyle{definition}
\newtheorem{defn}{Definition}

\newtheorem{example}{Example}
\newtheorem{problem}{Problem}

\theoremstyle{remark}
\newtheorem{rem}{Remark}

\begin{document} 
\title{Combinatorics of pedigrees}  
\author{Bhalchandra D. Thatte} 
\address{
Bio\-mathematics Research Centre \\
Mathematics and Computer Science Building \\
University of Canter\-bury \\
Private Bag 4800 \\
Christ\-church, New Zealand
} 
\email{bdthatte@gmail.com}
\keywords{pedigrees, reconstruction, enumeration}
\subjclass[2000]{Primary: 05C60}
\date{10 September, 2006} 
\begin{abstract}
	A pedigree is a directed graph in which each vertex 
	(except the founder vertices) has two parents. 
	The main result in this paper is a construction of
	an infinite family of counter examples to a reconstruction
	problem on pedigrees, thus negatively answering a question
	of Steel and Hein. Some positive reconstruction results
	are also presented.
	The problem of counting distinct (mutually non-isomorphic)
	pedigrees is considered. The known lower and upper bounds 
	on the number of pedigrees are improved upon, and their 
	relevance to pedigree reconstruction from DNA sequence data is 
	discussed. It is shown that the information theoretic bound
	on the number of segregating sites in the sequence data
	that is minimally essential for reconstructing pedigrees
	would not significantly change with improved enumerative
	estimates.
\end{abstract}
\maketitle

\section{Introduction}
\label{sec-intro}
A {\em general pedigree} $\mathcal{T}(X_0)$ 
on a set $X_0$, is a finite directed graph on 
a vertex set $V$ that satisfies the following conditions:
\begin{enumerate}
	\item each vertex has out-degree 0 or 2;
	\item $X_0$ is a subset of $V$, and each vertex in $X_0 $ has in-degree 0;
	\item there are no isolated vertices.
\end{enumerate}

The vertices with out-degree 0 are called the {\em founders}.
The vertices in $X_0$ are called {\em extant}. 
The cardinality of $X_0$ is called the {\em order} of the pedigree. 
Note that $X_0$ is a subset of the set of vertices with in-degree 0.

A {\em discrete generation pedigree} on $X_0$ is a pedigree on vertex 
set $V=\cup_{i=0}^d X_i$, where $X_i$ are disjoint sets, $X_d$ is the set 
of founders, and every vertex $u$ in $X_i; i < d$ has outgoing arcs $uv$ 
and $uw$ to vertices $v$ and $w$, respectively, in $X_{i+1}$. 
In this case, $d$ is the {\em depth} of the pedigree.

If there is an arc from a vertex $u$ to a vertex $v$, then $v$ is called 
a parent of $u$, and $u$ is called a child of $v$. If there is a directed 
path from a vertex $u$ to a vertex $v$ in a pedigree, then $v$ is said to be 
an {\em ancestor} of $u$, and $u$ is said to be a {\em descendent} of 
$v$. Trivially, each vertex is its own ancestor as well as its own descendent,
but not its own parent or child. If there is a directed path $u-u_1-\ldots 
u_k$ then $u_k$ is called a $k$-th grandparent of $u$, and $u$ is called 
a $k$-th grandchild of $u_k$.

A pedigree $\mathcal{P}(X_0)$ with vertex set $U$ is said to be {\em 
isomorphic} to a pedigree $\mathcal{Q}(Y_0)$ with vertex set $V$ if 
there is a one-one map $f:U\rightarrow V$ such that $u_1-u_2$ is an arc 
in $\mathcal{P}(X_0)$ if and only if $f(u_1)-f(u_2)$ is an arc in 
$\mathcal{Q}(Y_0)$. Although this is a standard definition of graph 
isomorphism, we will be interested in pedigrees in which the extant 
vertices are labelled. Therefore, if $X_0 = Y_0 = \{x_i; 1 \leq i \leq 
n\}$ then we will be interested only in isomorphisms $\pi$ for which 
$\pi(x_i) = x_i$ for all $1 \leq i \leq n$.

A motivation to study pedigrees comes from biology, where one is
interested in reconstructing pedigrees of populations. But it is hoped
that the main result in this paper - the non-reconstructibility of
pedigrees from sub-pedigrees - will also be of interst to combinatorialists
interested in the well known reconstruction conjectures.

Steel and Hein \cite{sh2006} posed and partially solved reconstruction 
and enumeration questions about pedigrees. Motivated by results in 
phylogenetics, a natural question to ask is: is a pedigree determined 
up to isomorphism from the pairwise distances between extant vertices? A 
pair of extant vertices $x$ and $y$ in a pedigree may have several 
common ancestors, therefore, it is assumed that all possible distances 
(in the undirected sense) between all pairs of extant vertices are 
given. Such a question is not expected to have a positive answer, as 
demonstrated by a counter example in \cite{sh2006}. Despite the counter 
example, variations of this question are definitely significant in 
evolutionary biology. Steel and Hein considered the following weaker question.

Let $\mathcal{P}(X_0)$ be a pedigree. A {\em sub-pedigree} 
$\mathcal{P}(Y)$ of $\mathcal{P}(X_0)$ is obtained by deleting every 
vertex in $\mathcal{P}(X_0)$ that has no descendent in $Y$. Now if 
sub-pedigrees on all two-element subsets of $X_0$ are given up to 
isomorphism, can we construct the sub-pedigree on $X_0$ up to 
isomorphism? Steel and Hein presented a counter example in their paper. 
They posed the following problem.

\begin{problem}
\label{prob-recon} 

Is there an integer $r>2$ such that every pedigree $\mathcal{P}(X_0)$ of 
order $n > r$ determined up to isomorphism if all its sub-pedigrees 
$\mathcal{P}(Y)$ such that $|Y| = r$ are given up to isomorphism?

\end{problem}

Combinatorialists familiar with the reconstruction conjectures might be 
tempted to dismiss this question, therefore, it must be pointed out that 
the set $X_0$ in a pedigree is labelled. In other words, ``a 
sub-pedigree $\mathcal{P}(Y)$ given up to isomorphism'' is to be 
interpreted as a pedigree in which all vertices except the ones in $Y$ 
are unlabelled. The following definitions are introduced to make this 
remark more formal.

\begin{defn}

Let $ n > r > 2 $ be positive integers. Let $\mathcal{T}(X_0)$ and 
$\mathcal{U}(X_0)$ be two pedigrees of order $n$. The two pedigrees are 
said to be $r$-{\em hypomorphic} to each other if for every $Y \subset 
X_0; |Y| = r$, there is an isomorphism $\pi_Y$ from the sub-pedigree 
$\mathcal{T}(Y)$ of $\mathcal{T}(X_0)$ to the sub-pedigree 
$\mathcal{U}(Y)$ of $\mathcal{U}(X_0)$ such that $\pi_Y(x) = x$ for all $x 
\in X_0$. A pedigree $\mathcal{T}(X_0)$ is said to be $r$-{\em 
reconstructible} if for every pedigree $\mathcal{U}(X_0)$ that is 
$r$-hypomorphic to $\mathcal{T}(X_0)$, there is an isomorphism $\pi$ from 
$\mathcal{T}(X_0)$ to $\mathcal{U}(X_0)$ such that $\pi(x) = x$ for all $x 
\in X_0$.

\end{defn}

\begin{problem}
\label{prob-recon1}

Is there an integer $r > 2$ such that all pedigrees of order $n > r$ are 
$r$-reconstructible?

\end{problem}

In Section~\ref{sec-recon}, we present a family of counter examples as 
well as a few positive results on constant population size pedigrees. We 
prove that for every $n>3$, there are pedigrees of order $n$ that are 
not even $(n-1)$-reconstructible. The problem of classification of 
non-reconstructible pedigrees remains open, and we suspect that it might 
have an algebraic structure similar to the Nash-Williams' lemma in edge 
reconstruction theory, see \cite{nash-williams-1978}.

Steel and Hein considered the question of enumerating mutually 
non-isomorphic pedigrees of a fixed depth. A lower bound on the number 
of distinct pedigrees implies, by an information theoretic argument, a 
lower bound on the number of segregating DNA sites that would be 
necessary in order to reconstruct the pedigree of a population from the 
sequence data. In Section~\ref{sec-enum} we prove tighter 
lower and upper bounds, and show that the information theoretic lower 
bound does not increase much. Steel and Hein leave the problem of 
enumerating general pedigrees open. Here we enumerate general pedigrees 
as well, and again show that purely information theoretic arguments as
in their paper are not sufficient to show that general pedigrees would
necessarily require significantly more segregating sites for their
reconstruction from the sequence data.

\section{Reconstruction of pedigrees.}
\label{sec-recon}
\subsection{A negative result}
\label{subsec-negative}
We solve Problem ~\ref{prob-recon} negatively by constructing an
infinite family of pairs of non-isomorphic pedigrees 
that have correspondingly isomorphic sub-pedigrees. That is, we prove 
the following

\begin{thm}
\label{thm-counter}
For every $n > 2$, there are non-isomorphic pedigrees $\mathcal{T}(X_0)$ 
and $\mathcal{U}(X_0)$ of order $n$ that are $(n-1)$-hypomorphic.
\end{thm}

\begin{proof}

The proof is divided in two cases. The case $n=3$ gives the basic idea, 
which is then generalised to arbitrary values of $n$.
	
\noindent {\bf Case $n = 3$}. 

Consider the non-isomorphic graphs $K_{1,3} $ and $K_3$. Let the edges 
of both graphs be arbitrarily labelled $e_1,e_2,e_3$, where, following 
the standard graph theoretic convention, an edge is a set of two 
vertices. It is clear that $K_{1,3} - e_i \cong K_3 - e_i$ for all $i$, 
where $-e_i$ denotes deletion of the edge $e_i$ and the resulting 
isolated vertices. Now suppose that the end vertices of each edge $e_i$ 
are parents of the vertex $x_i \in X_0$ in each of the pedigrees 
$\mathcal{T}(X_0)$ and $\mathcal{U}(X_0)$. Then the pedigrees 
$\mathcal{T}(\{x_i,x_j\})$ and $\mathcal{U}(\{x_i,x_j\})$ are isomorphic 
for all $i,j$, but the pedigrees $\mathcal{T}(X_0)$ and 
$\mathcal{U}(X_0)$ are not isomorphic. This example proves the theorem 
for $n = 3$. The pedigrees $\mathcal{T}(X_0)$, $\mathcal{U}(X_0)$,
and their sub-pedigrees are shown in Figure~\ref{fig-base-case}.
\begin{figure}[ht]
\includegraphics{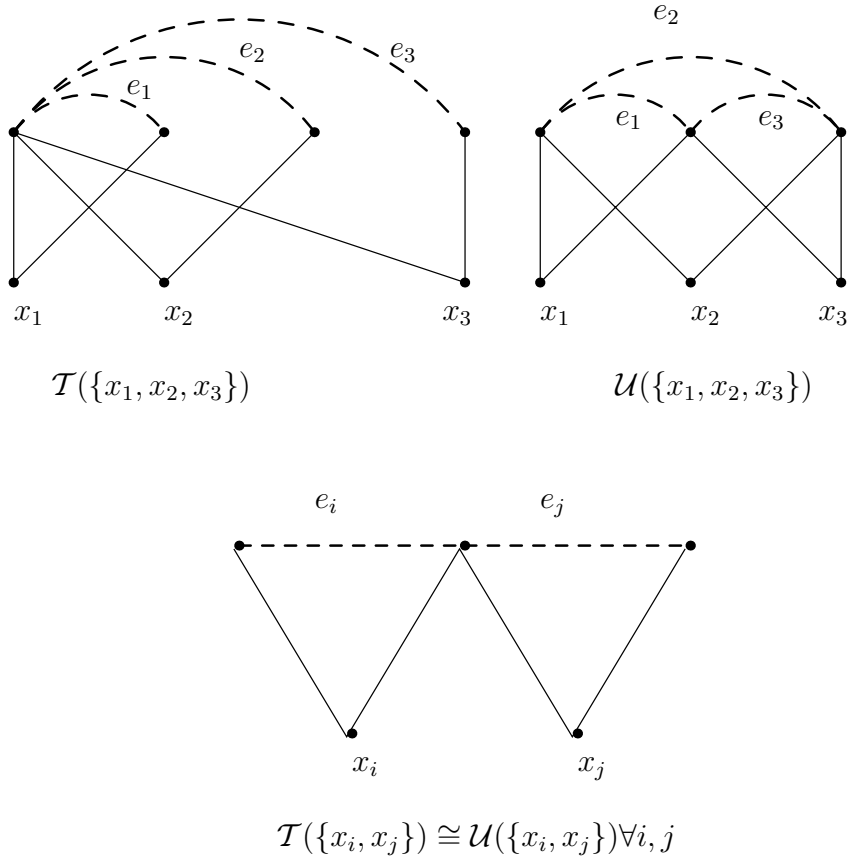}
\vspace{0.2in}
\caption[]{Pedigrees based on $K_{1,3}$ and $K_3$}
\label{fig-base-case}
\end{figure}

\noindent {\bf Case $n > 3$} \\
We have to construct hypergraphs that play the role that $K_{1,3}$ and 
$K_3$ play above. We construct a hypergraph $G$ with edge set $\{g_i; 1 \leq 
i \leq n\}$ and a hypergraph $H$ with edge set $\{h_i; 1 \leq i \leq n\}$
such that the following conditions are satisfied.

\begin{enumerate} 
\item $G \ncong H$ 
\item For each $i; 1 \leq i \leq n$, $G-g_i \cong H-h_i$; moreover, there is 
an isomorphism between $G-g_i$ and $H-h_i$ that {\em preserves the edge 
order}, that is, vertices in an edge $g_j$ in $G-g_i$ are mapped to 
vertices in $h_j$ under such an isomorphism, for each $j \neq i$.

\end{enumerate}

Once such hypergraphs are constructed, we treat each edge in each 
hypergraph as a founder set. We construct pedigrees $T_i$ on founder 
sets $g_i$, and pedigrees $U_i$ on founder sets $h_i$ such that

\begin{enumerate}

\item pedigrees $T_i; 1 \leq i \leq n$ are vertex-disjoint except 
possibly for their founder sets $g_i$;

\item pedigrees $U_i; 1 \leq i \leq n$ are vertex-disjoint except 
possibly for their founder sets $h_i$;

\item pedigrees $T_i$ and $U_i$ are correspondingly isomorphic; 
moreover, for all $i,j; i \neq j$, an isomorphism between $G-g_j$ and 
$H-h_j$ that preserves the edge order extends to an isomorphism between 
$T_i$ and $U_i$;

\item each of the pedigrees $T_i$ and $U_i$ contains
exactly one extant vertex $x_i$.

\end{enumerate}

The resulting pedigrees $\cup_{i=1}^n T_i$ and $\cup_{i=1}^n U_i$ are 
non-isomorphic (since the hypergraphs $G$ and $H$ are non-isomorphic) 
but their sub-pedigrees are correspondingly isomorphic.

\noindent {\bf Construction of hypergraphs $G$ and $H$}

The required hypergraphs are constructed by a simple application of 
linear algebra.

Let each integer in $\{0,2^{n}-1\}$ be written in base 2 as an $n$-digit 
number by padding sufficiently many zeros on the left. We count its 
digits from the right. The set of $n$-digit binary numbers is denoted by 
$[2^n]$. The $i$'th digit of a number $k$ is denoted by $k(i)$, and the 
number obtained by setting the $i$'th digit of $k$ to 0 (or 1) is 
denoted by $k(i\leftarrow 0)$ (or, respectively, $k(i\leftarrow 1)$). 
The number of ones and the number of zeros in $k$ are denoted by 
$\#1(k)$ and $\#0(k)$, respectively.

The isomorphism class of a hypergraph $G$ with edge set $\{g_i; 1 \leq i 
\leq n\}$ may be represented by a list of integers $a(k); k \in [2^n]$, 
where $a(k)$ is the number of vertices in $\cap_{i=1}^n f_i$, where $f_i 
= g_i$ if $k(i) = 1$, and $f_i = \bar{g_i}$, (that is, the complement of 
$g_i$), if $k(i) = 0$. In other words, we have to only specify the 
number of vertices in each region of the Venn diagram of $g_i; 1\leq i 
\leq n$. Let the list of integers $b(k); k \in [2^n]$ similarly denote 
the isomorphism class of $H$.

The condition $G-g_i \cong H-h_i$ for $1 \leq i \leq n$, with an 
isomorphism between them that preserves the edge order, may be expressed 
as
\begin{equation} 
	\label{eq-cube}
a(k(i\leftarrow 0))+a(k(i\leftarrow 1)) 
= b(k(i\leftarrow 0)) + b(k(i\leftarrow 1)); k \in [2^n]
\end{equation}

Since we are interested in non-isomorphic hypergraphs $G$ and $H$, we 
must find solutions to the above equations so that $a(k) \neq b(k)$ for 
some $k \in [2^n]$.

We verify that
\begin{align}
a(k) & = 1, &  b(k) = 0 && \text{when $k$ has even number of 1's} \\
a(k) & = 0, &  b(k) = 1 && \text{when $k$ has odd number of 1's} \notag 
\end{align}
satisfy Equations~(\ref{eq-cube}).

It can be easily verified that $K_{1,3}$ and $K_3 \cup K_1$ do in fact
satisfy the above solutions, where we include an isolated vertex in one
of the graphs purely for algebraic convenience.

Now on we write $[2^n] = [2^n]_e \cup [2^n]_o$, where $[2^n]_e$ is 
the set of integers having an even number of 1's in their binary
representation, and $[2^n]_o$ is the set of integers having 
an odd number of 1's in their binary representation.
In this notation, the hypergraphs $G$ and $H$ are described as follows:
the set $[2^n]_e$ is the set of vertices of $G$, and a vertex 
$k\in [2^n]_e$ is in $g_i$ if and only if $k(i) = 1$. 
Similarly, the set $[2^n]_o$ is the vertex set of $H$, and a vertex
$k\in [2^n]_o$ is in edge $h_i$ if and only if $k(i) = 1$.

The vertex $k=0$ is in $G$, but is an isolated vertex, and is included at
this stage only for algebraic convenience, and may be deleted after completing
the construction of non-reconstructible pedigrees.

It is clear that $G$ and $H$ are non-isomorphic, since each of them has 
$2^{n-1}$ vertices, but $G$ has the isolated vertex 0, while $H$ has no 
isolated vertex. What is an isomorphism between $G-g_i$ and $H-h_i$? An 
edge order preserving isomorphism from $G-g_i$ to $H-h_i$ must map 
vertices in a region of the Venn diagram of $\cup_{j|j\neq i}g_j$ to the 
corresponding region of the Venn diagram of $\cup_{j|j\neq i}h_j$.
Consider any $k \in [2^n]$. The vertex $k(i\leftarrow 0)$ is in 
$g_j$ for some $j \neq i$ if and only if the vertex $k(i\leftarrow 1)$ is
in $h_j$, because the two vertices differ only in their $i$'th digit.
Therefore, if $k(i\leftarrow 0)$ is in $G$, then an edge order
preserving isomorphism between $G-g_i$ and $H-h_i$ must map the vertex
$k(i\leftarrow 0)$ to the vertex $k(i\leftarrow 1)$. 
Similarly, if the vertex $k(i\leftarrow 1)$ is in $G$, then an edge 
order preserving isomorphism between $G-g_i$ 
and $H-h_i$ must map the vertex $k(i\leftarrow 1)$ to the vertex $k(i 
\leftarrow 0)$. Moreover, this isomorphism is unique. On the standard 
hypercube on $[2^n]$, each vertex in $G-g_i$ is mapped to its neighbour 
along the $i$'th axis, which is in $H-h_i$.

\begin{example}
Let $n=4$, and let the hypergraphs $G$ and $H$ be defined on vertex
sets $[2^4]_e$ and $[2^4]_o$ as follows:
\begin{align}
g_1 = \{0011, 0101, 1001, 1111\}, \,\,&  
g_2 = \{0011, 0110, 1010, 1111\},\notag \\
g_3 = \{0101, 0110, 1100, 1111\}, \,\,&
g_4 = \{1001, 1010, 1100, 1111\},\notag \\
h_1 = \{0001, 0111, 1011, 1101\}, \,\,& 
h_2 = \{0010, 0111, 1011, 1110\},\notag \\
h_3 = \{0100, 0111, 1101, 1110\}, \,\,&
h_4 = \{1000, 1011, 1101, 1110\},\notag 
\end{align}
where $g_i$ are the edges of $G$ and $h_i$ are the edges of $H$.
The isomorphism $\pi_1$ from $G-g_1$ to $H-h_1$ that preserves the edge order is
given by $\pi_1(0000) = 0001$, $\pi_1(0011) = 0010$, $\pi_1(0101) = 0100$,
$\pi_1(1001) = 1000$, $\pi_1(0110) = 0111$, $\pi_1(1010) = 1011$, $\pi_1(1100) = 1101$,
$\pi_1(1111) = 1110$. Observe that $\pi_1(g_2) = h_2$, $\pi_1(g_3) = h_3$, and
$\pi_1(g_4) = h_4$ under this map.
\end{example}

\noindent {\bf Construction of $T_i$ and $U_i$ }

As stated earlier, for each $i$, pedigrees $T_i$ and $U_i$ must be so 
constructed that (the unique) edge order preserving isomorphism between 
$G-g_j$ and $H-h_j$ extends to an isomorphism between $T_i$ and $U_i$ 
for all $j \neq i$.

Let a balanced binary tree $T_i$ be defined so that $x_i$ is its root 
and $g_i$ is its set of leaves. By convention, the root $x_i$ is the 
lowest vertex (at depth 0) in $T_i$, and the leaves are the highest 
vertices (at depth $n-2$) in $T_i$. For a vertex $t$ in $T_i$, let 
$T_i(t)$ be the subtree of $T_i$ induced by $t$ and all vertices in 
$T_i$ that are above $t$. Let $t_0$ and $t_1$ be the parents of $t$. The 
subtree $T_i(t)$ is a union of subtrees $L(t)$ and $R(t)$, where $L(t)$ 
is induced by vertices $t$, $t_0$, and all vertices above $t_0$, and the 
subtree $R(t)$ is induced by vertices $t$, $t_1$, and all vertices above 
$t_1$. We call $L(t)$ the {\em left subtree} at $t$, and $R(t)$ the {\em 
right subtree} at $t$.

Let $i_1, i_2,\ldots , i_{n-1}$ be the integers $1 \leq j \leq n; j \neq i$
in arbitrary order. The vertices in $g_i$ are grouped in such a 
way that for each vertex $t$ at depth $k; 0 \leq k \leq n-3$, if a 
vertex $p\in g_i$ is a leaf of $L(t)$ then $p(i_{k+1})=0$, and if a 
vertex $p\in g_i$ is a leaf of $R(t)$ then $p(i_{k+1})=1$.

The vertices in $h_i$ are partitioned, and a binary tree $U_i$ is 
constructed analogously for the same ordering $i_j; 1 \leq j \leq n-1$.

For $n=5$ and $i=5$ and the ordering $i_1 = 2, i_2 = 3, i_3 = 1, i_4 = 
4$, the trees $T_5$ and $U_5$ are shown in Figure ~\ref{fig-general-case}.
\begin{figure}[ht]
\includegraphics[scale=0.92]{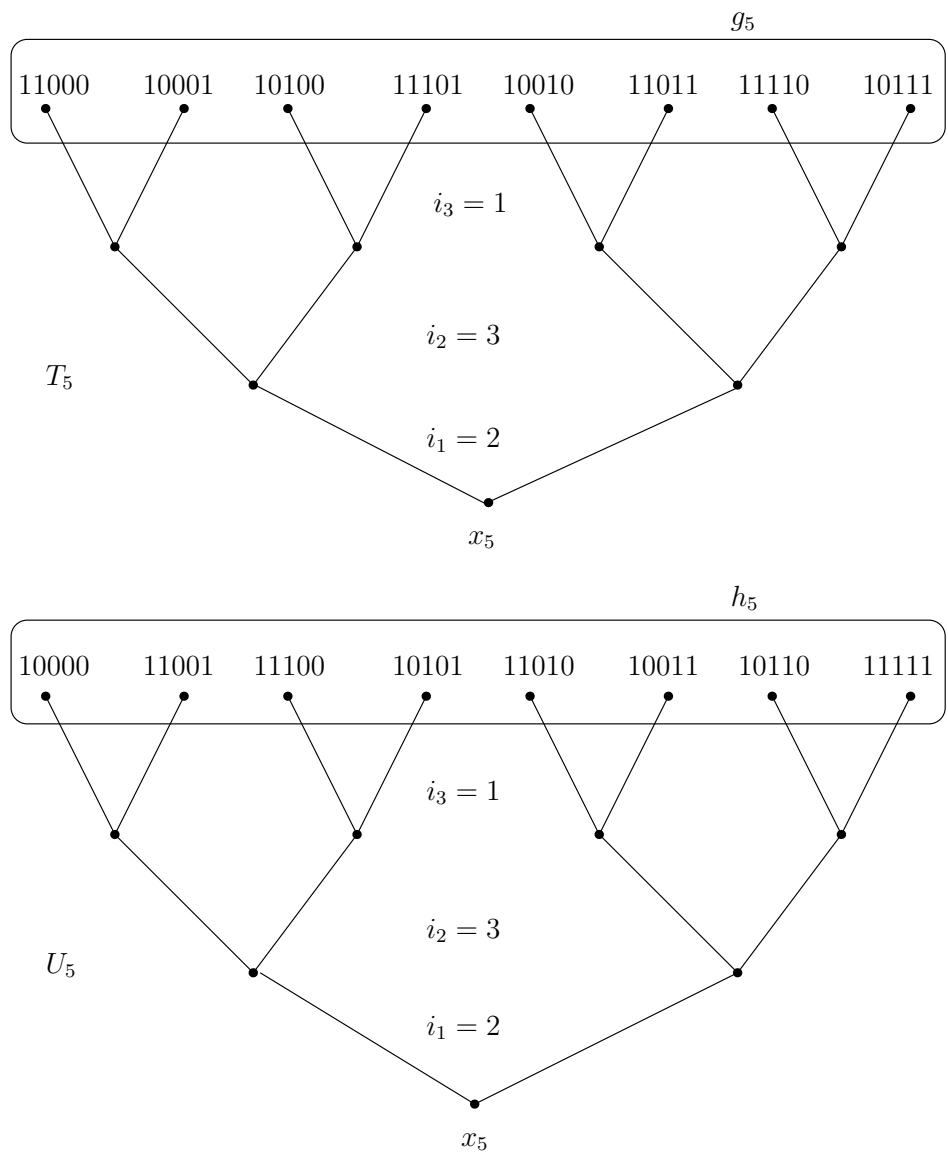}
\vspace{0.2in}
\caption[]{Binary pedigrees $T_5$ and $U_5$ for $n=5$ }
\label{fig-general-case}
\end{figure}

We show that for every $j \neq i$, the unique isomorphism between 
$G-g_j$ and $H-h_j$ extends to an isomorphism between $T_i$ and $U_i$.

Let $\bar{b}=(b_1, \ldots , b_j); 0 \leq j \leq n-2$ be a $j$-tuple of 
0's and 1's. Extending a notation introduced earlier, let $g_i(\bar{b})$ 
denote the set $\{k\in g_i| k(i_1) = b_1, k(i_2) = b_2, \dots , k(i_j) = 
b_j\}$, which is the set of leaves of a binary subtree of $T_i$ rooted 
at the vertex $t(\bar{b})$ at depth $j$. For example, when $\bar{b} = 
(0)$, $g_i(\bar{b})$ is the set of leaves above the left parent of 
$x_i$, and $t(\bar{b})$ is the left parent of $x_i$. The set 
$h_i(\bar{b})$ and the vertex $u(\bar{b})$ are analogously defined for 
$U_i$. By convention, an empty tuple $\bar{b}$ defines the sets $g_i$ 
and $h_i$, and the trees $T_i$ and $U_i$, rooted at $x_i$; and a tuple 
$\bar{b}$ of length $n-2$ defines singleton subsets $\{t(\bar{b})\}$ of 
$g_i$, and $\{u(\bar{b})\}$ of $h_i$. A tuple $\bar{b}$ of length $n-2$ also 
uniquely determines the digits $t(\bar{b})(i_{n-1})$ and 
$u(\bar{b})(i_{n-1})$, since we know that $\#1(t(\bar{b}))$ is even and 
$\#1(u(\bar{b}))$ is odd. Also, if $t(\bar{b}(i_{n-1})) = 1$ then 
$u(\bar{b}(i_{n-1})) = 0$, and if $t(\bar{b}(i_{n-1})) = 0$ then 
$u(\bar{b}(i_{n-1})) = 1$. Therefore, the map $t(\bar{b}) 
\longrightarrow u(\bar{b})$ for all tuples of length at most $n-2$ 
extends the isomorphism between $G-g_{i_{n-1}}$ and $H-h_{i_{n-1}}$.

Let $\bar{b}$ be a tuple as above. Extending the notation $k(i\leftarrow 
0)$ to $\bar{b}$, we define $\bar{b}(i\leftarrow 0)$ to be the tuple 
obtained by setting $b_i=0$ in $\bar{b}$, and $\bar{b}(i\leftarrow 1)$ 
to be the tuple obtained by setting $b_i=1$ in $\bar{b}$.

Let $\bar{b}$ be a tuple of length $n-2$ and $j \leq n-2$. By an 
argument as in the above paragraph, we have

\begin{enumerate}

\item if $t(\bar{b}(j\leftarrow 0)(i_{n-1})) = 1$ then 
$u(\bar{b}(j\leftarrow 1)(i_{n-1})) = 1$;

\item if $t(\bar{b}(j\leftarrow 0)(i_{n-1})) = 0$ then 
$u(\bar{b}(j\leftarrow 1)(i_{n-1})) = 0$;

\item if $t(\bar{b}(j\leftarrow 1)(i_{n-1})) = 1$ then 
$u(\bar{b}(j\leftarrow 0)(i_{n-1})) = 1$;

\item if $t(\bar{b}(j\leftarrow 1)(i_{n-1})) = 0$ then 
$u(\bar{b}(j\leftarrow 0)(i_{n-1})) = 0$.

\end{enumerate}

Therefore, for each $j; 1\leq j \leq n-2$, the map defined by 

\begin{enumerate}

\item $t(\bar{b}) \longrightarrow u(\bar{b})$ for all tuples of length 
at most $j-1$;

\item $t(\bar{b}(j\leftarrow 0))\longrightarrow u(\bar{b}(j\leftarrow 1))$ 
for all tuples of length at least $j$; and

\item $t(\bar{b}(j\leftarrow 1))\longrightarrow u(\bar{b}(j\leftarrow 0))$ 
for all tuples of length at least $j$

\end{enumerate} extends the isomorphism between $G-g_j$ and $H-h_j$. 
Observe that this map sends the vertices in the left subtree 
$L(t(\bar{b}))$ in $T_i$ to the vertices in the right subtree 
$R(u(\bar{b}))$ in $U_i$, and the vertices in the right subtree 
$R(t(\bar{b}))$ in $T_i$ to the vertices in the left subtree 
$L(u(\bar{b}))$ in $U_i$, for each tuple $\bar{b}$ of length $j$.
Since the trees $T_i$ and $T_j$ (and trees $U_i$ and $U_j$) are disjoint
except for their founders for all $i \neq j$, the isomorphism
between $G-g_j$ and $H-h_j$ extends to an isomorphism between
pedigrees $\mathcal{T}(X_0\backslash\{x_j\})$ and
$\mathcal{U}(X_0\backslash\{x_j\})$ for all $j$.

\begin{figure}[ht]
\includegraphics[scale=0.88]{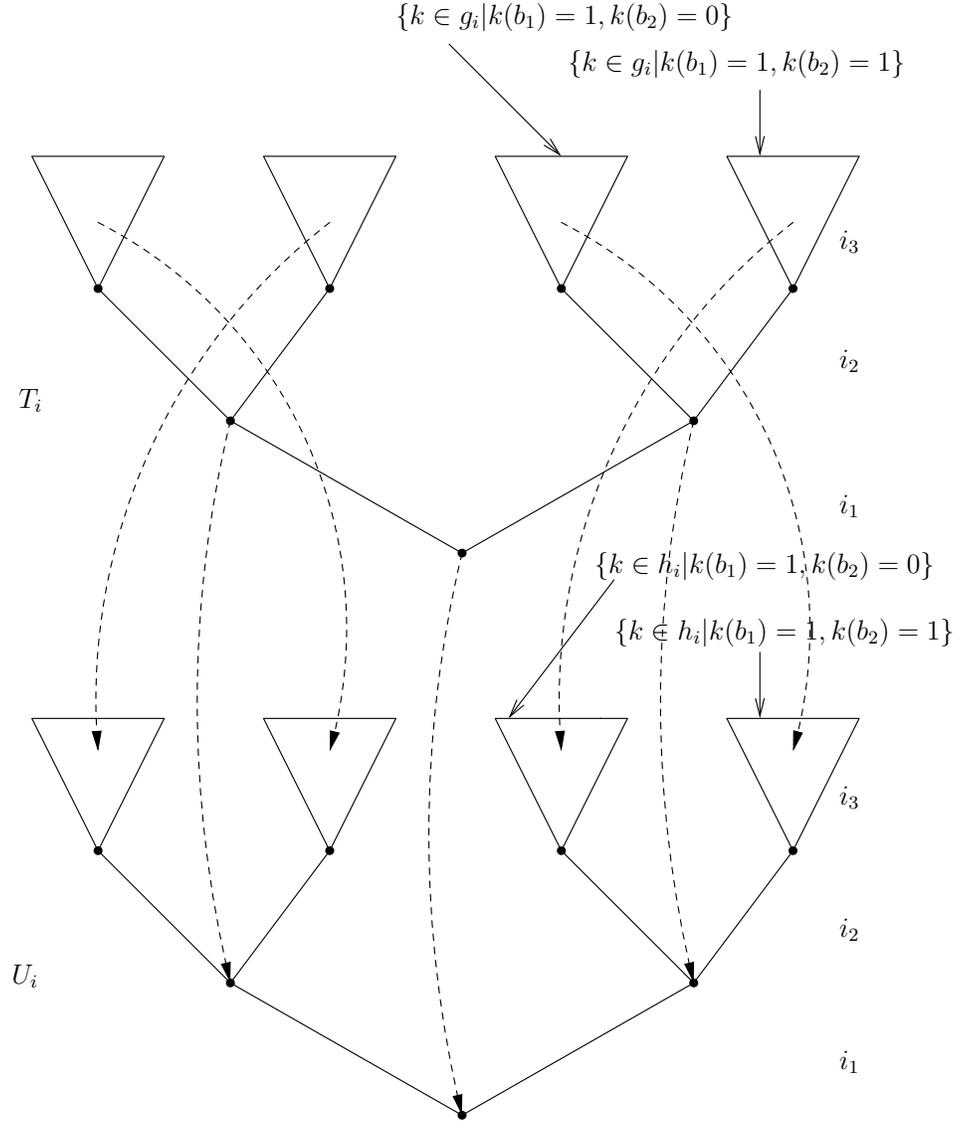}
\vspace{0.2in}
\caption[]{Isomorphism between $G-g_{i_2}$ and $H-h_{i_2}$ extends to an isomorphism between $T_i$ and $U_i$ }
\label{fig-isomorphism}
\end{figure}

An isomorphism between $G-g_{i_2}$ and $H-h_{i_2}$ that
extends to an isomorphism between $T_i$ and $U_i$; $i \neq i_2$
is schematically shown in Figure~\ref{fig-isomorphism}. 
\end{proof}

\begin{figure}[ht]
\includegraphics[scale=0.87]{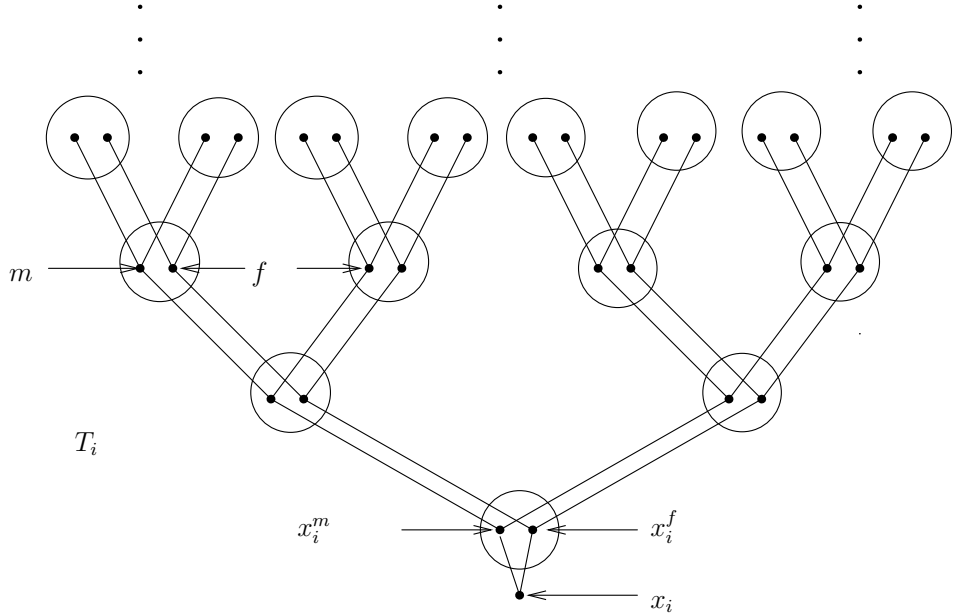}
\vspace{0.2in}
\caption[]{Construction of a pedigree with a valid gender labelling }
\label{fig-duplicate}
\end{figure}

\begin{rem}
	The pedigrees constructed above do not admit a valid gender labelling.
	That is, we cannot assign labels $m$ (male) and $f$ (female) to
	all vertices so that each vertex (except founders) has one male
	parent and one female parent. For example, in the $n=3$ case,
	$K_3$ is not a bipartite graph, so a valid gender labelling is
	impossible. But the examples can be easily modified to create
	non-reconstructible pedigrees that also admit valid gender
	labels. Each vertex in a pedigree constructed above may be duplicated,
	and one vertex may be treated male and the other female, as
	shown in Figure~\ref{fig-duplicate}. At the bottom of the tree
	$T_i$ (or $U_i$), the vertex $x_i$ is duplicated as vertices
	$x_i^m$ and $x_i^f$, and the new vertex $x_i$ is a child of
	$x_i^m$ and $x_i^f$.  
\end{rem}

\begin{rem}
	\label{rem-minimal}
	Let $k, k^\prime \in [2^n]$ be any two adjacent vertices on the
	hypercube. From Equation~(\ref{eq-cube}), if $a(k) - b(k) =  p > 0$
	for some $p$, then $b(k^\prime) - a(k^\prime) = p$, regardless
	of which digit $k$ and $k^\prime $ differ at. In fact, by
	connectivity of the hypercube, we have $a(r) -b(r) = p$ for all
	vertices $r \in [2^n]$ that are at even distance from $k$ on the
	hypercube, and $b(r) - a(r) = p$ for all vertices $r \in [2^n]$ that are
	at odd distance from $k$ on the hypercube. This further implies
	that the hypergraphs $G$ and $H$ constructed in the above counter example
	have a special structure: for each $i ; 1 \leq i \leq n$, $|g_i|
	= |h_i| \geq 2^{n-2}$. Let $G(d)$ be the hypergraph with edge set
	$\{g_i(d); 1 \leq i \leq n\}$, where $g_i(d)$ is the set of
	grandparents of $x_i$ at depth $d$ in the pedigree
	$\mathcal{T}(X_0)$, and let $H(d)$ and $h_i(d)$ be similarly defined
	for the pedigree $\mathcal{U}(X_0)$, then the hypergraphs $G(d)$
	and $H(d)$ must be isomorphic whenever $d < n-2$.
\end{rem}

We end this subsection with a conjecture motivated by the observations made
in Remark~\ref{rem-minimal}.
\begin{conjecture}
	\label{con-minimal}
	The counter example constructed in Theorem~\ref{thm-counter} is
	minimal.  In other words, if a pedigree $\mathcal{T}(X_0)$ of order
	$n$ is not $(n-1)$-reconstructible then it has depth at least $n-2$, and
	there are at least $2^{n-1}$ ancestors at depth $n-2$. 
\end{conjecture}

\begin{rem}
	\label{rem-nash}
	Let $G$ and $H$ be simple graphs with edge sets $E(G) = \{g_i; 1\leq i
	\leq m\}$ and $E(H) = \{h_i; 1 \leq i \leq m\}$, respectively, 
	such that $G-g_i\cong H-h_i$ for all $1 \leq i \leq m$. Then the
	edge reconstruction conjecture states that $G\cong H$ provided
	$m > 3$. The condition $m > 3$ is required since $K_{1,3}$ and
	$K_3$ - the graphs used as the base case of our construction of
	non-reconstructible pedigrees - are not {\em edge reconstructible}.
	Although no counter examples are yet known, Nash-Williams
	\cite{nash-williams-1978} proved a characterisation of
	(hypothetical) counter examples to edge reconstruction. 
	His characterisation was based on a generalisation of ideas
	earlier introduced by Lov\'asz \cite{lovasz-1972}.
	Without going into details, we note that the counter examples
	presented here have certain similarities with the
	characterisation by Nash-Williams. It may be possible to
	exploit such similarities to prove Conjecture~\ref{con-minimal}. 
\end{rem}

\subsection{A positive result}
\label{subsec-positive}

Let $\mathcal{T}(X_0)$ be a discrete generation pedigree on $X_0$
of order $n > 2$.
Let $S_{n-1}(\mathcal{T})= \{\mathcal{T}(Y)| Y \subset X_0, \, |Y| = n-1\}$.
Consider the edge labelled (multi) graph $G_1$ whose vertex 
set is $X_1$ (that is, the vertices at depth $1$), and vertices $x,y\in 
X_1$ are joined by an edge $e_i$ if they are the parents of $x_i$.

\begin{lem}
	\label{lem-same-parents}
	If there are vertices $x_i$ and $x_j$ in $X_0$ that have the same
	parents, then $\mathcal{T}(X_0)$ is uniquely determined by
	$S_{n-1}(\mathcal{T})$.
\end{lem}

\begin{proof}
The situation in the lemma is recognised by looking at 
$\mathcal{T}(X_0\backslash x_k)$, where $x_k\not \in \{x_i , x_j\}$. Now
$\mathcal{T}(X_0)$ is uniquely obtained from $\mathcal{T}(X_0\backslash x_i)$
by joining $x_i$ to the parents of $x_j$.
\end{proof}

\begin{lem}
	\label{lem-cycle}
	If $n > 3$ and if $G_1$ contains a cycle then $\mathcal{T}(X_0)$ is 
	uniquely determined $S_{n-1}(\mathcal{T})$.
\end{lem}

\begin{proof}

Let $e_i$ be an edge in a cycle in $G_1$. The end vertices of $e_i$ are
the two parents of $x_i$. Since the set of half brothers 
of $x_i$ is known from the collection $S_{n-1}$, the parents of $x_i$ 
are uniquely recognised in $\mathcal{T}(X\backslash x_i)$. Note that we need the 
condition $n> 3$ because otherwise we would get a counter example based 
on $G_1 \cong K_3$ or $G_1 \cong K_{1,3}$.

\end{proof}

\begin{cor}
	If $|X_1| \leq n$ and $ n > 3$ then $S_{n-1}(\mathcal{T})$ determines 
	$\mathcal{T}(X_0)$ up to congruence.
\end{cor}
\begin{proof}
	If no two vertices in $X_0$ have the same two parents then $G_1$
	has $n$ simple edges (that is no two edges are parallel edges),
	and there is a cycle in $G_1$. 
\end{proof}

We end this section with another conjecture.
\begin{conjecture}
	\label{con-logn}
	Discrete generation pedigrees of order $n$ that have a constant
	population in each generation are $r$-reconstructible for
	$ r > \log n$.
\end{conjecture}
This conjecture is true if Conjecture~\ref{con-minimal} is true.
For suppose that Conjecture~\ref{con-minimal} is true but
Conjecture~\ref{con-logn} is not true, and that there is a pedigree of
order $n$ that is not $r$-reconstructible for some $r > \log n$. 
Therefore, for some $r > \log n$, there is a sub-pedigree of order $r+1$
that is not $r$-reconstructible.
Such a sub-pedigree must have depth at least $r-1$, and must have at least
$2^{r}$ vertices at depth $r-1$, implying that $r \leq \log n$.
Thus if $r > \log n $ then we have a contradiction, therefore,
all sub-pedigrees of order $r+1$ are $r$-reconstructible when $r > 
\log n $, and we can complete the reconstruction inductively.

\section{Enumeration of pedigrees}
\label{sec-enum} 

Let $N(n,d)$ be the number of distinct (mutually non-isomorphic)
discrete generation pedigrees of depth $d$ with $n$ vertices in
each generation. As before, the extant vertices are assumed to be
labelled, and other vertices are assumed to be unlabelled. 

In a general pedigree, the depth of a vertex $u$ is the largest integer
$k$ for which $u$ is a $k$'th grandparent of an extant vertex.
The depth of a pedigree is the largest integer $d$ for which
there is a vertex of depth $d$ in the pedigree.
Let the number of distinct general pedigrees of depth $d$
with constant number $n$ of vertices at each depth be $M(n,d)$.  

The purpose of this section is to derive lower and upper bounds on
$N(n,d)$ and $M(n,d)$. The bounds are relevant to an information
theoretic argument that was used by Steel and Hein in the context
of a reconstruction question. 

\begin{thm}
	\label{thm-bounds}
\begin{align}
   	 \left(\frac{(n-1)n^{n-2}}{2}\right)^d 
	 &\leq & N(n,d) & \leq & \binom{n}{2}^{nd}
	 \label{eq-bound1} \\
	 \frac{(n-1)n^{n-2}}{2}\prod_{k=0}^{d-2}((n/2)(d-1-k))^n
	 & \leq & M(n,d) & \leq & \binom{nd-1}{2}^{nd}
	 \label{eq-bound2}
\end{align}
\end{thm}

\begin{proof}

Let $\mathcal{P}(X_0)$ be a discrete generation pedigree of depth
$d$ on $X_0$.  Let $X_i$ be the set of vertices at depth $i$.
Let $|X_i| = n$ for all $i; 0 \leq i \leq d$.
For each $i; 1\leq i \leq d$, define a graph $G_i$ as follows:
the vertex set of $G_i$ is $X_i$, and $\{u,v\}$ is an edge in
$G_i$ if $u$ and $v$ have a child in $X_{i-1}$.
Thus $1 \leq e(G_i) \leq n$ for $1 \leq i \leq d$,
where $e(G)$ denotes the number of edges of a graph $G$.
We restrict ourselves to bipartite graphs $G_i$ so that it is
possible to assign valid gender labels to the vertices
of pedigrees.  

Let $S(n,k)$ denote the Sterling number of the second kind.
There are $S(n,k)$ partitions of $X_0$ in groups of siblings, 
where {\em siblings} are vertices that share both parents.
If the vertices of $G_1$ are labelled then there are $k!$ ways
of assigning the groups of siblings to pairs of parents.
Therefore, each labelled graph $G_1$ gives $S(n,k)k!$ labelled
pedigrees of depth 1.
Some of those pedigrees may be isomorphic to each other since there
may be automorphisms of $G_1$ that permute the edges of $G_1$ non-trivially.
Therefore, the number of distinct pedigrees of depth 1 that can be
obtained from a labelled graph $G_1$ is given by
\[
N(n,d,G_1) = \frac{S(n,k)k!}{|\text{aut}LG_1|},
\]
where $LG_1$ denotes the line graph of $G_1$, and  $\text{aut}G$ denotes
the automorphism group of a graph $G$.
If every non-trivial automorphism of $G_1$ permutes the edges of $G_1$
non-trivially then the number of distinct pedigrees of depth 1 
that can be obtained from $G_1$ is given by
\[
N(n,d,G_1) = \frac{S(n,k)k!}{|\text{aut}G_1|}.
\]
Each non-trivial automorphism of a graph $G$ permutes the edges of $G$
non-trivially if and only if $G$ has no isolated edges and not more
than one isolated vertices. Therefore,
\[
N(n,1) \geq \sum_G \frac{S(n,e(G))e(G)!}{|\text{aut}G|},
\]
where the summation is over all distinct bipartite graphs $G$
having $n$ vertices, at least 1 and at most $n$ edges, at most
one isolated vertex, and no isolated edges.

Pedigrees of depth 1 considered above have the additional
property that they have no non-trivial automorphisms that fix
each vertex in $X_0$, implying that the vertices of $X_1$ are
distinguishable in such pedigrees.
Therefore,
\[
N(n,d) \geq \left(\sum_G \frac{S(n,e(G))e(G)!}{|\text{aut}G|}\right)^d,
\]
where the summation is over all graphs of the type described above.

Summing over only graphs that have $n-1$ edges, we have
\[
\frac{S(n,e(G))e(G)!}{|\text{aut}G|} = \frac{\binom{n}{2}(n-1)!}{|\text{aut}G|}
\]
But $n!/|\text{aut}G|$ is the number of labelled graphs isomorphic to
$G$. Therefore, summing over trees, we get
\[
N(n,d) \geq \left(\frac{(n-1)n^{n-2}}{2}\right)^d.
\]

The upper bound on $N(n,d)$ is obtained by counting fully labelled
pedigrees that do not even possibly admit a valid gender labelling.

We derive a lower bound on $M(n,d)$ by enumerating a
special subclass of general pedigrees that is described next.
Consider pedigrees of depth $d$ and order $n$ that satisfy the conditions:
\begin{enumerate}
	\item there are $n$ vertices at each depth $k\leq d$, 
	\item each vertex at depth $k\leq d-2$ has exactly one 
		parent at depth $k+1$,
	\item distinct vertices at depth $k \leq d-2$ have distinct
		parents at depth $k+1$, 
	\item at each depth $k; k \leq d-1$, there are $n/2$ vertices
		of each gender,
	\item the pedigree of depth 1 induced by vertices in
		$X_{d-1}\cup X_d$ has no non-trivial automorphisms
		that fix vertices in $X_{d-1}$. 
\end{enumerate}

The conditions imply that given any vertex $v$ at depth $k; k\leq d-1$
there is a unique path of length $k$ beginning at some vertex $u$ in
$X_0$ and ending at $v$.
Therefore, vertices at depth at most $d-1$ are distinguishable.
The last condition above makes the vertices at depth $d$ distinguishable
as well. Therefore, no two pedigrees described by the above conditions are
isomorphic.  This allows us to derive a lower bound on $M(n,d)$.
\[
M(n,d) \geq \frac{(n-1)n^{n-2}}{2}\prod_{k=0}^{d-2}((n/2)(d-1-k))^n,
\]
where the first factor is a lower bound on the number of distinct
pedigrees of depth 1 that are induced by $X_{d-1}\cup X_d$, vertices in
$X_{d-1}$ being labelled. For a vertex at depth $k\leq d-2$, the parent
that is not at depth $k+1$ may be chosen from the $(n/2)(d-1-k)$
distinguishable vertices at depth $k+2$ or more. This explains the
second factor. 

An upper bound on $M(n,d)$ is obtained by counting the number
of labelled directed graphs in which each vertex has out-degree 2.
\end{proof}

\begin{rem}
	Steel and Hein give the information theoretic argument
	that if there are $s$ segregating sites in DNA sequences
	obtained from $n$ extant individuals, then there are $4^{ns}$
	possible combinations of sequences. Therefore, $4^{ns}$
	must be at least $N(n,d)$ (or $M(n,d)$) depending on what
	assumptions are made about pedigrees) to be able to reconstruct
	their pedigree up to depth $d$.
	They derive a lower bound on $s$ given by $(d/3)\log n$ for
	reconstruction of discrete generation constant population size pedigrees.
	They comment that in reality the number of sites required is
	likely to be much higher due to under-counting  of isomorphism
	classes and due to the stochastic nature of sequence evolution.
	Theorem~\ref{thm-bounds} gives an information theoretic lower 
	bound on $s$ that is about $(d/2)\log n$ for
	discrete generation constant population size pedigrees, and a bound
	of about $(d/2)\log (nd)$ for general pedigrees. Moreover, the
	bounds based on the upper bounds on $N(n,d)$ and $M(n,d)$ are
	only about $d\log n$ and $d\log (nd)$, respectively, for
	discrete generation and general pedigrees.  
\end{rem}

\begin{rem}
If we assume that no vertex at depth $k$ has a parent
at depth more than $k+t+1$ then we have
\[
M(n,d) \geq \frac{(n-1)n^{n-2}}{2}
\prod_{k=0}^{d-t-1}(nt/2)^n\prod_{k=d-t}^{d-2}(n(d-k-1)/2)^n \\
\]
This gives a lower bound of about $(d/2) \log(nt)$ on the
number of segregating sites required for pedigree reconstruction.
\end{rem}
\subsection* { Acknowledgements } I would like to thank Mike Steel for many
useful conversations, and for hosting me at the University of
Canterbury. Oliver Will and Mike Steel read the manuscript carefully and
made many helpful comments. I am supported by the Allan Wilson Centre for
Molecular Ecology and Evolution, New Zealand.

\bibliographystyle{plain}

\begin{thebibliography}{1}

\bibitem{lovasz-1972}
L.~Lov{\'a}sz.
\newblock A note on the line reconstruction problem.
\newblock {\em J. Combinatorial Theory Ser. B}, 13:309--310, 1972.

\bibitem{nash-williams-1978}
C.~St. J.~A. Nash-Williams.
\newblock The reconstruction problem.
\newblock In Lowell~W. Beineke and Robin~J. Wilson, editors, {\em Selected
  topics in graph theory}, pages 205--236. Academic Press Inc. [Harcourt Brace
  Jovanovich Publishers], London, 1978.

\bibitem{sh2006}
M.~Steel and J.~Hein.
\newblock Reconstructing pedigrees: a combinatorial perspective.
\newblock {\em Journal of Theoretical Biology}, 2006.

\end{thebibliography}

\end{document}